\documentclass[12pt,a4paper,twoside]{article}

\newcommand{\pageformat}[6]{\setlength{\hoffset}{-1in}
                  \setlength{\voffset}{-1in}
                  \addtolength{\hoffset}{#5}
                            \addtolength{\voffset}{#6}
                            \setlength{\oddsidemargin}{#1}
                            \setlength{\evensidemargin}{#2}
                            \setlength{\textwidth}{\paperwidth}
                  \addtolength{\textwidth}{-\oddsidemargin}
                  \addtolength{\textwidth}{-\evensidemargin}
                  \addtolength{\textwidth}{-\marginparsep}
                  \addtolength{\textwidth}{-\marginparwidth}
                            \setlength{\topmargin}{#3}
                            \setlength{\textheight}{\paperheight}
                  \addtolength{\textheight}{-\topmargin}
                  \addtolength{\textheight}{-\headheight}
                  \addtolength{\textheight}{-\headsep}
                  \addtolength{\textheight}{-\footskip}
                  \addtolength{\textheight}{-#4}}
\pageformat{2cm}{3cm}{25mm}{25mm}{1pt}{0pt}

\usepackage{ifthen}
\newboolean{article}
    \setboolean{article}{true}
\newboolean{report}
\newboolean{book}
\newboolean{letter}
\newboolean{german}
\newboolean{italian}
\newboolean{nobaselinestretch}
\newboolean{nosectionappendix}
\newboolean{oldtoc}
\newboolean{nosectionequn}
\newboolean{notheorem}

\ifthenelse{\boolean{german}}{
    \usepackage{german}}{}

\usepackage[latin1]{inputenc}

\usepackage{amsmath}
\usepackage{amssymb}
\usepackage[mathscr]{eucal}

\ifthenelse{\boolean{notheorem}}{}{
    \usepackage{theorem}}

\ifthenelse{\boolean{nobaselinestretch}}{}{
    \renewcommand{\baselinestretch}{1.25}}

\newenvironment{env}[2]{\begin{#1}#2\end{#1}}{}
    \newcommand{\beq}[1]{\begin{env}{equation}{#1}}
    \newcommand{\beqn}[1]{\begin{env}{equation*}{#1}}
    \newcommand{\bal}[1]{\begin{env}{align}{#1}}
    \newcommand{\baln}[1]{\begin{env}{align*}{#1}}
    \newcommand{\bga}[1]{\begin{env}{gather}{#1}}
    \newcommand{\bgan}[1]{\begin{env}{gather*}{#1}}
    \newcommand{\bflal}[1]{\begin{env}{flalign}{#1}}
    \newcommand{\bflaln}[1]{\begin{env}{flalign*}{#1}}
    \newcommand{\bmu}[1]{\begin{env}{multline}{#1}}
    \newcommand{\bmun}[1]{\begin{env}{multline*}{#1}}
    \newcommand{\bsp}[1]{\begin{env}{split}{#1}}

    \newcommand{\eeq}{\end{env}}
    \newcommand{\eeqn}{\end{env}}
    \newcommand{\eal}{\end{env}}
    \newcommand{\ealn}{\end{env}}
    \newcommand{\ega}{\end{env}}
    \newcommand{\egan}{\end{env}}
    \newcommand{\eflal}{\end{env}}
    \newcommand{\eflaln}{\end{env}}
    \newcommand{\emu}{\end{env}}
    \newcommand{\emun}{\end{env}}
    \newcommand{\esp}{\end{env}}

\newcommand{\lf}{\vspace{2ex}}

\renewcommand{\bf}[1]{\textbf{#1}}
\renewcommand{\it}[1]{\textit{#1}}

\renewcommand{\sf}[1]{\textsf{#1}}

\renewcommand{\tt}[1]{\texttt{#1}}
\newcommand{\hl}[1]{\bf{\it{#1}}}

\newcommand{\msf}[1]{\text{\small$\sf{#1}$}}

\newcommand{\eus}[1]{\mathscr{#1}}
\newcommand{\euf}[1]{\mathfrak{#1}}
\newcommand{\bb}[1]{\mathbb{#1}}

\newcommand{\nbd}[1]{$#1$\nobreakdash--}

\newcommand{\wt}[1]{\widetilde{#1}}

\newcommand{\vt}{\vartheta}

\newcommand{\om}{\omega}

\newcommand{\norm}[1]{\left\lVert#1\right\rVert}

\newcommand{\Babs}[1]{\Bigl\lvert#1\Bigr\rvert}
\newcommand{\Bnorm}[1]{\Bigl\lVert#1\Bigr\rVert}
\newcommand{\snorm}[1]{\norm{\smash{#1}}}

\newcommand{\family}[1]{\left(#1\right)}

\newcommand{\bfam}[1]{\bigl(#1\bigr)}
\newcommand{\Bfam}[1]{\Bigl(#1\Bigr)}
\newcommand{\AB}[1]{\langle#1\rangle}
\newcommand{\bAB}[1]{\bigl\langle#1\bigr\rangle}
\newcommand{\BAB}[1]{\Bigl\langle#1\Bigr\rangle}
\newcommand{\CB}[1]{\{#1\}}
\newcommand{\bCB}[1]{\bigl\{#1\bigr\}}
\newcommand{\BCB}[1]{\Bigl\{#1\Bigr\}}
\newcommand{\SB}[1]{[#1]}
\newcommand{\bSB}[1]{\bigl[#1\bigr]}

\newcommand{\LO}[1]{(#1]}

\newcommand{\set}[2][]{
    \ifthenelse{\equal{#1}{}}{
        \CB{#2}}{
        \CB{#1~|~#2}}}
\newcommand{\bset}[2][]{
    \ifthenelse{\equal{#1}{}}{
        \bCB{#2}}{
        \bCB{#1~|~#2}}}
\newcommand{\Bset}[2][]{
    \ifthenelse{\equal{#1}{}}{
        \BCB{#2}}{
        \BCB{#1~\big|~#2}}}
\newcommand{\zero}{\CB{0}}

\DeclareMathOperator{\id}{\normalfont\msf{id}}

\newcommand{\C}{\bb{C}}

\newcommand{\N}{\bb{N}}

\newcommand{\R}{\bb{R}}

\newcommand{\sB}{\eus{B}}

\newcommand{\eH}{\euf{H}}

\newcommand{\I}{{I\!\!\!\;I}}

\ifthenelse{\boolean{nosectionequn}}{}{
    \numberwithin{equation}{section}
    }

\ifthenelse{\boolean{article}\or\boolean{letter}\or\boolean{nosectionequn}}{
    \setboolean{nosectionappendix}{true}}{}
\ifthenelse{\boolean{nosectionappendix}}{}{
    \renewcommand{\appendix}{
        \chapter*{\appendixname}
        \addcontentsline{toc}{chapter}{\appendixname}
        \renewcommand{\thesection}{\Alph{section}}
        \setcounter{section}{0}}}
   
\ifthenelse{\boolean{report}\or\boolean{book}}{
    }{}

\ifthenelse{\boolean{notheorem}}{}{
        \newcommand{\mnname}{Mathematical note.}
        \newcommand{\enname}{End of the note.}
        \newcommand{\definame}{Definition.}
        \newcommand{\propname}{Proposition.}
        \newcommand{\lemname}{Lemma.}
        \newcommand{\exname}{Example.}
        \newcommand{\exername}{Exercise.}
        \newcommand{\remname}{Remark.}
        \newcommand{\obname}{Observation.}
        \newcommand{\thmname}{Theorem.}
        \newcommand{\corname}{Corollary.}
        \newcommand{\proofname}{Proof.}
    \ifthenelse{\boolean{german}}{
        \renewcommand{\mnname}{Mathematische Notiz.}
        \renewcommand{\enname}{Ende der Notiz.}
        \renewcommand{\exname}{Beispiel.}
        \renewcommand{\exername}{Übung.}
        \renewcommand{\remname}{Bemerkung.}
        \renewcommand{\obname}{Beobachtung.}
        \renewcommand{\thmname}{Satz.}
        \renewcommand{\corname}{Korollar.}
        \renewcommand{\proofname}{Beweis.}}{}
    \ifthenelse{\boolean{italian}}{
        \renewcommand{\mnname}{Nota matematica.}
        \renewcommand{\enname}{Fina della nota.}
        \renewcommand{\definame}{Definizione.}
        \renewcommand{\propname}{Proposizione.}
        \renewcommand{\exname}{Esempio.}
        \renewcommand{\exername}{Esercizio.}
        \renewcommand{\remname}{Nota.}
        \renewcommand{\obname}{Osservazione.}
        \renewcommand{\thmname}{Teorema.}
        \renewcommand{\corname}{Corollario.}
        \renewcommand{\proofname}{Dimostrazione.}

       \renewcommand{\appendixname}{Appendice}

       }{}
    \theoremheaderfont{\normalfont\bfseries}
    \theoremstyle{change}
        \theorembodyfont{\rmfamily}
            \newtheorem{emp}{}[section]
                \newcommand{\bemp}[1][]{
                    \begin{emp}\hskip-\labelsep\bf{#1}\hskip\labelsep}
                \newcommand{\eemp}{\end{emp}}
\newtheorem{itemp}[emp]{}
                \newcommand{\bitemp}[1][]{
                    \begin{itemp}\hskip-\labelsep\bf{#1}\hskip\labelsep\normalfont\itshape}
                \newcommand{\eitemp}{\end{itemp}}
            \newtheorem{mn}[emp]{\mnname}
                \newcommand{\bnm}{\begin{mn}~\begin{quotation}\renewcommand{\baselinestretch}{1}\small\noindent\ignorespaces}
                \newcommand{\enm}{\end{quotation}\hfill\bf{\enname}\end{mn}}
            \newtheorem{ex}[emp]{\exname}
                \newcommand{\bex}{\begin{ex}}
                \newcommand{\eex}{\end{ex}}
            \newtheorem{exer}[emp]{\exername}
                \newcommand{\bexer}{\begin{exer}}
                \newcommand{\eexer}{\end{exer}}
            \newtheorem{defi}[emp]{\definame}
                \newcommand{\bdefi}{\begin{defi}}
                \newcommand{\edefi}{\end{defi}}
            \newtheorem{rem}[emp]{\remname}
                \newcommand{\brem}{\begin{rem}}
                \newcommand{\erem}{\end{rem}}
            \newtheorem{ob}[emp]{\obname}
                \newcommand{\bob}{\begin{ob}}
                \newcommand{\eob}{\end{ob}}
        \theorembodyfont{\normalfont\itshape}
            \newtheorem{thm}[emp]{\thmname}
                \newcommand{\bthm}{\begin{thm}}
                \newcommand{\ethm}{\end{thm}}
            \newtheorem{prop}[emp]{\propname}
                \newcommand{\bprop}{\begin{prop}}
                \newcommand{\eprop}{\end{prop}}
            \newtheorem{cor}[emp]{\corname}
                \newcommand{\bcor}{\begin{cor}}
                \newcommand{\ecor}{\end{cor}}
            \newtheorem{lem}[emp]{\lemname}
                \newcommand{\blem}{\begin{lem}}
                \newcommand{\elem}{\end{lem}}
\newenvironment{empn}[1]{\lf\noindent\bf{#1}\ignorespaces\hskip\labelsep}{\lf}
		\newcommand{\bempn}[1]{\begin{empn}{#1}}
		\newcommand{\eempn}{\end{empn}}
		\newcommand{\bitempn}[1]{\begin{empn}{#1}\normalfont\itshape}
		\newcommand{\eitempn}{\end{empn}}
                \newcommand{\bnmn}{\begin{empn}{\mnname}~\begin{quotation}\renewcommand{\baselinestretch}{1}\small\noindent\ignorespaces}
                \newcommand{\enmn}{\end{quotation}\hfill\bf{\enname}\end{empn}}
		\newcommand{\bexn}{\begin{empn}{\exname}}
		\newcommand{\eexn}{\end{empn}}
		\newcommand{\bexern}{\begin{empn}{\exername}}
		\newcommand{\eexern}{\end{empn}}
		\newcommand{\bdefin}{\begin{empn}{\definame}}
		\newcommand{\edefin}{\end{empn}}
		\newcommand{\bremn}{\begin{empn}{\remname}}
		\newcommand{\eremn}{\end{empn}}
		\newcommand{\bobn}{\begin{empn}{\obname}}
		\newcommand{\eobn}{\end{empn}}
		\newcommand{\bthmn}{\bitempn{\thmname}}
		\newcommand{\ethmn}{\eitempn}
		\newcommand{\bpropn}{\bitempn{\propname}}
		\newcommand{\epropn}{\eitempn}
		\newcommand{\bcorn}{\bitempn{\corname}}
		\newcommand{\ecorn}{\eitempn}
		\newcommand{\blemn}{\bitempn{\lemname}}
		\newcommand{\elemn}{\eitempn}
\newcommand{\qedsymbol}{~\rule[-0.35mm]{2mm}{2mm}}
    \newcounter{proof}[emp]
    \newenvironment{Proof}[1]{
        \vspace{1ex}
        \renewcommand{\item}[1][\stepcounter{proof}(\roman{proof})]%
            {##1\hskip\labelsep}
        \noindent\textsc{#1\hskip\labelsep}}{
        \nolinebreak\qedsymbol}
    \newcommand{\proof}[1][\proofname]{
        \begin{Proof}{#1}\ignorespaces}
    \newcommand{\qed}{\end{Proof}}
    \newcommand{\noqed}{
        \renewcommand{\qedsymbol}{}
        \end{Proof}}}
    \ifthenelse{\boolean{italian}}{
        \renewcommand{\proofname}{Dimostrazione.}}{}

\usepackage[varg]{txfonts}

\renewcommand{\thefootnote}{[\arabic{footnote}]}

\usepackage[hypertex]{hyperref}

\begin{document}

\bibliographystyle{amsalpha}

\title{Constructing Proper Markov Semigroups\\for Arveson Systems\thanks{This work has been supported by research funds of the Dipartimento S.E.G.e S.\ of University of Molise and of the Italian MIUR (PRIN 2007).}}

\author{Michael Skeide}

\date{June 2010}

\maketitle

\begin{abstract}
\noindent
{\renewcommand{\thefootnote}{}}We show that the Markov semigroup obtained by Floricel in \cite{Flo08} compressing the \nbd{E_0}semi\-group of Skeide \cite{Ske06}, does not consist of endomorphisms. It, therefore, cannot be the tail flow of an \nbd{E_0}semigroup. As a corollary of our result, Floricel's construction will allow to get examples of proper type III Markov semigroups that are not tensor products of simpler ones, provided we find type III Arveson systems that do not factor into tensor products.
\end{abstract}

\noindent
Algebraically, an \hl{Arveson system} is a family $E^\otimes=\bfam{E_t}_{t\in(0,\infty)}$ of infinite-dimensional separable Hilbert spaces $E_t$ with unitaries $u_{t,s}\colon E_t\otimes E_s\rightarrow  E_{t+s}$ such that the product $(x_t,y_s)\mapsto x_ty_s:=u_{t,s}(x_t\otimes y_s)$ is associative. Technically, the bundle $\bfam{E_t}_{t\in(0,\infty)}$ is required to be a Borel bundle isomorphic to the trivial Borel bundle $(0,\infty)\times\eH$ for some Hilbert space $\eH$, and the product is required measurable.

Suppose we find a Hilbert space $K(\ne\zero)$ and an (again measurable) family of unitaries $w_t\colon E_t\otimes K\rightarrow K$ such that the product $(x_t,y)\mapsto x_ty:=w_t(x_t\otimes y)$ iterates associatively with the product of the Arveson system. Then the maps $\vt_t\colon a\mapsto w_t(\id_t\otimes a)w_t^*$ on $\sB(K)$ (plus $\vt_0:=\id_{\sB(K)}$) form an \hl{\nbd{E_0}semigroup} $\vt=\bfam{\vt_t}_{t\in\R_+}$ (that is, a \nbd{\sigma}weakly continuous semigroup of normal unital endomorphisms). Such a family $\bfam{w_t}_{t\in(0,\infty)}$ has been called a \hl{right dilation} of $E^\otimes$ in Skeide \cite{Ske06} and a \hl{unitary resolution} in Floricel \cite{Flo08}. It is equivalent to the notion of \hl{nondegenerate} or \hl{essential} representation of an Arveson system. Arveson \cite{Arv89} associated with every \nbd{E_0}semigroup an Arveson system, and it is not difficult to show that the Arveson system of $\vt$ is $E^\otimes$.

Arveson also proved in \cite{Arv90} that every Arveson system admits a right dilation. Skeide \cite{Ske06} presented an elementary direct construction of a right dilation, and Floricel \cite{Flo08} generalized that further. The first ingredient of the construction in \cite{Ske06}, is a right dilation $\breve{w}_n\colon E_n\otimes\breve{K}\rightarrow\breve{K}$ $(n\in\N)$ of the \hl{discrete subsystem} $\bfam{E_n}_{n\in\N}$ of $E^\otimes$. Such a right dilation can be obtained from any unit vector $\om_1\in E_1$ as an inductive limit $\breve{K}$ over $E_n$ with respect to the inductive system $E_n\rightarrow E_n\om_1^m\subset E_{n+m}$. It is not difficult to check that the factorization $E_n\otimes E_m\rightarrow E_{n+m}$ survives the ``limit'' $m\to\infty$, giving $\breve{w}_n$. Moreover, all $\om_1^n\in E_n$ end up in the same unit vector $\breve{\om}\in\breve{K}$, which fulfills $\om_1^n\breve{\om}=\breve{\om}$. One may check that $\AB{\breve{\om},\bullet\breve{\om}}$ is an absorbing invariant vector state for the discrete \nbd{E_0}semigroup $\breve{\vt}$ on $\sB(\breve{K})$. In particular, the projections $\breve{\vt}_n(\breve{\om}\breve{\om}^*)$ increase to the identity. See \cite[Appendix]{Arv89} or \cite[Section 5]{BhSk00} for details.

We rest a moment to clarify some details about tensor products with direct integrals and operations on them. If $H_1,H_2$ are separable Hilbert spaces, then for $0\le a<b<\infty$ we will understand by
\beqn{
\int_a^b(H_1\otimes E_\alpha\otimes H_2)\,d\alpha
}\eeqn
the family of measurable, square integrable sections $X=\bfam{x_\alpha}_{\alpha\in\LO{a,b}}$ with $x_\alpha\in H_1\otimes E_\alpha\otimes H_2$. Since $\bfam{E_\alpha}_{\alpha\in\LO{a,b}}$ is Borel isomorphic to $\LO{a,b}\times\eH$, it is clear how this has to be interpreted. In particular, $\int_a^b(H_1\otimes E_\alpha\otimes H_2)\,d\alpha\cong L^2(\LO{a,b},H_1\otimes \eH\otimes H_2)$ by the Borel isomorphism. It is clear that
\beqn{
H_1\otimes\Bfam{\int_a^bE_\alpha\,d\alpha}\otimes H_2
~\cong~
\int_a^b(H_1\otimes E_\alpha\otimes H_2)\,d\alpha
}\eeqn
via $x_1\otimes\bfam{y_\alpha}_{\alpha\in\LO{a,b}}\otimes x_2\mapsto\bfam{x_1\otimes y_\alpha\otimes x_2}_{\alpha\in\LO{a,b}}$, because
\beqn{
L^2(\LO{a,b},H_1\otimes\eH\otimes H_1)
~\cong~
L^2\LO{a,b}\otimes H_1\otimes\eH\otimes H_2
~\cong~
H_1\otimes L^2(\LO{a,b},\eH)\otimes H_2.
}\eeqn
Recall that for $x_2\in H_2$ the operator $\id_1\otimes x_2^*\in\sB(H_1\otimes H_2,H_1)$ is defined by setting $(\id_1\otimes x_2^*)(y_1\otimes y_2)=y_1\AB{x_2,y_2}$. It is the adjoint of the operator $\id_1\otimes x_2\colon y_1\mapsto y_1\otimes x_2$.

\bpropn
Let $X=\bfam{x_\alpha}_{\alpha\in\LO{a,b}}\in\int_a^bE_\alpha\,d\alpha$. Then the operator $\id_1\otimes X^*$ acts on sections $Y=\bfam{y_\alpha}_{\alpha\in\LO{a,b}}\in\int_a^b(H_1\otimes E_\alpha)\,d\alpha$ as
\beqn{
(\id_K\otimes X^*)Y
~=~
\int_a^b(\id_1\otimes x_\alpha^*)y_\alpha\,d\alpha
}\eeqn
in the sense of Bochner integral of \nbd{H_1}valued functions.
\epropn

Similar statements are true for $\int_a^b(E_\alpha\otimes H_2)\,d\alpha$ and $\int_a^b(H_1\otimes E_\alpha\otimes H_2)\,d\alpha$.

\bcorn
$(\id_1\otimes X'X^*)Y=\family{~\Bfam{\int_a^b(\id_1\otimes x_\beta^*)y_\beta\,d\beta}~\otimes~x'_\alpha}_{\alpha\in\LO{a,b}}$.
\ecorn

\proof[Proof of the proposition.~]
Evaluate the operator on the dense set of elements of the form
\beqn{
Y
~=~
\sum_{i=1}^n\bfam{h_i\otimes(\I_{\LO{a_i,b_i}}(\alpha)y'_\alpha)}_{\alpha\in\LO{a,b}}
}\eeqn
($h_i\in H_1$, $Y'\in\int_a^bE_\alpha\,d\alpha$, and the $\LO{a_i,b_i}$ forming a partition of $\LO{a,b}$) and extend it in \nbd{L^2}norm.\qed

\lf
Note that $\breve{K}$ and the right dilation $\breve{w}_n$ of the discrete subsystem of $E^\otimes$ to $\breve{K}$ coincide with $\wt{K}_u$ and $\wt{W}_{u,n}$ in \cite{Flo08}, if one puts the $s>0$ in \cite[Section 3.1]{Flo08} equal to $s=1$ and $u\in E_s$ equal to $u=\om_1$. The vector $\breve{\om}$, in \cite{Flo08} is denoted by $\SB{u}$.

Put $K:=\Bfam{\int_0^1E_\alpha\,d\alpha}\otimes\breve{K}=\int_0^1(E_\alpha\otimes\breve{K})\,d\alpha$. Both \cite{Ske06} and \cite{Flo08} define right dilations $w_t$ and $W_{u,t}$ ($s=1$ and $u=\om_1$ as before), respectively, of $E^\otimes$ to $K$. We do not know, if the two right dilations coincide, or if the \nbd{E_0}semigroups $\vt$ and $\rho$, respectively, determined by them coincide. However, we know that they coincide for integer $t=n$ and this is enough for our purposes.

Indeed, for $t=n\in\N$ the right dilation $w_t$ defined in \cite[Equation (3.1)]{Ske06}%
\footnote[1]{\label{1}
Well, actually in \cite[Equation (3.1)]{Ske06} a left dilation is defined. By \cite[Theorem 3.3]{Ske06}, reversing the orders in all tensor products one gets a right dilation. This is, how \cite[Equation (3.1)]{Ske06} must be interpreted here.
}%
 ~acts as
\beqn{
w_n
\colon
x_n\otimes \bfam{y_\alpha\otimes\breve{z}}_{\alpha\in\LO{0,1}}
~=~
\bfam{x_n\otimes y_\alpha\otimes\breve{z}}_{\alpha\in\LO{0,1}}
~\longmapsto~
\Bfam{\,(\id_\alpha\otimes\breve{w}_n)\,\bSB{(u_{\alpha,n}^*(x_ny_\alpha)\otimes\breve{z})}\,}_{\alpha\in\LO{0,1}}.
}\eeqn
And this is precisely what the definition of $W_{u,t,l}$ in \cite[Equation (3.18)]{Flo08} according to the equation between Equations (3.16) and (3.17) in \cite{Flo08} gives for $W_{u,n}:=W_{\om_1,n,0}$.

Recall that if $\vt$ is an \nbd{E_0}semigroup on $\sB(K)$ and if $Q\in\sB(K)$ is an \hl{increasing projection} (that is, if $\vt_t(Q)\ge Q$ for all $t$), then the maps $T_t\colon QaQ\mapsto Q\vt_t(QaQ)Q=Q\vt_t(a)Q$ form a \hl{Markov semigroup} $T$, that is, a semigroup of normal unital CP-maps, on $Q\sB(K)Q$. We call $T$ the Markov semigroup obtained from $\vt$ by \hl{compression} with $Q$.

We abbreviate $L:=\int_0^1E_\alpha\,d\alpha$, so that $K=L\otimes\breve{K}$. Define $Q:=\id_L\otimes\breve{\om}\breve{\om}^*\in\sB(K)$. We confirm \cite[Proposition 4.2]{Flo08}:

\blemn
$Q$ is increasing for $\vt$.
\elemn

For integer times $t=n$, the proof will be evident from an intermediate step in the proof of the following result. For arbitrary $t$ we would have to repeat the full definition of $w_t$ from \cite{Ske06}, and for the following proof it does not matter if the maps $T_t$ form a Markov semigroup.

\bthmn
The Markov semigroup $T$ does not consist of endomorphisms.
\ethmn

\proof
Note that $Q\sB(K)Q=\sB(L)\otimes\breve{\om}\breve{\om}^*\cong\sB(L)$. A normal, unital ($T_n$ is Markov!) endomorphism of $\sB(L)$ takes non-zero projections to non-zero projections. We shall show that there exists a rank-one projection $a\in Q\sB(K)Q$ such that $T_1(a)$ is not a non-zero projection.

Fix a unit vector $X=\bfam{x_\alpha}_{\alpha\in\LO{0,1}}\in L$, and define the rank-one projection $a:=XX^*\otimes\breve{\om}\breve{\om}^*\in Q\sB(H)Q$. The norm of the positive operator $T_n(a)\in\sB(L)\otimes\breve{\om}\breve{\om}^*$ is the supremum of the matrix elements $\AB{(Y\otimes\breve{\om}),T_n(a)(Y\otimes\breve{\om})}$ over all unit vectors $Y=\bfam{y_\alpha}_{\alpha\in\LO{0,1}}\in L$. First, we observe that $Q(Y\otimes\breve{\om})=Y\otimes\breve{\om}$. Next, we compute
\bmun{
w_n^*(Y\otimes\breve{\om})
~=~
\bfam{(u_{n,\alpha}^*u_{\alpha,n}\otimes\id_{\breve{K}})(\id_\alpha\otimes\breve{w}_n^*)(y_\alpha\otimes\breve{\om})}_{\alpha\in\LO{0,1}}
\\
~=~
\bfam{(u_{n,\alpha}^*u_{\alpha,n}\otimes\id_{\breve{K}})(y_\alpha\otimes\breve{\om}_n\otimes\breve{\om})}_{\alpha\in\LO{0,1}}
~=~
\bfam{u_{n,\alpha}^*(y_\alpha\om_1^n)\otimes\breve{\om}}_{\alpha\in\LO{0,1}}
}\emun
Finally,
\bmun{\tag{$\dagger$}\label{D}
\bAB{(Y\otimes\breve{\om}),(Q\vt_n(a)Q)(Y\otimes\breve{\om})}
~=~
\bAB{w_n^*(Y\otimes\breve{\om}),(\id_n\otimes XX^*\otimes\breve{\om}\breve{\om}^*)w_n^*(Y\otimes\breve{\om})}
\\
~=~
\bAB{\bfam{u_{n,\alpha}^*(y_\alpha\om_1^n)}_{\alpha\in\LO{0,1}},(\id_n\otimes XX^*)\bfam{u_{n,\alpha}^*(y_\alpha\om_1^n)}_{\alpha\in\LO{0,1}}}
\\
~=~
\BAB{
\int_0^1(\id_n\otimes x_\beta^*)u_{n,\beta}^*(y_\beta\om_1^n)\,d\beta,
\int_0^1(\id_n\otimes x_\gamma^*)u_{n,\gamma}^*(y_\gamma\om_1^n)\,d\gamma}
\\
~=~
\Bnorm{
\int_0^1(\id_n\otimes x_\gamma^*)u_{n,\gamma}^*(y_\gamma\om_1^n)\,d\gamma}^2.
}\emun
(At this point, replacing in the first two lines $XX^*$ with $\id_L$, we see that, indeed, $\vt_n(Q)\ge Q$; that proves the preceding Lemma for integer times: $T_n$ is, indeed, Markov.)

We put $n=1$ and we shall find a unit vector $X$ such that \eqref{D} is not bigger than a constant $M^2<1$ no matter what unit vector $Y$ we choose. Note that there exists a unit vector $z_1\in E_1$ such that the square root of \eqref{D} is given by
\beqn{
\BAB{z_1,\int_0^1(\id_1\otimes x_\gamma^*)u_{1,\gamma}^*(y_\gamma\om_1)\,d\gamma}
~=~
\int_0^1\AB{z_1x_\gamma,y_\gamma\om_1}\,d\gamma.
}\eeqn
Choose a  \hl{measurable ONB} for $E^\otimes$. By this, we mean a family $\bfam{e^m}_{m\in\N}$ of measurable sections $e^m=\bfam{e^m_\alpha}_{\alpha\in\LO{0,1}}$ such that for each $\alpha$ the family $\bfam{e^m_\alpha}_{m\in\N}$ is an ONB for $E_\alpha$. (Such a measurable ONB exists, because $E^\otimes$ is isomorphic to a trivial bundle.) Then the vectors $f^m_{1-\alpha}:=(\id_{1-\alpha}\otimes{e^m_\alpha}^*)\om_1\in E_{1-\alpha}$, $\alpha\in(0,1)$ depend measurably on $\alpha$ (all Hilbert spaces are separable) and fulfill $\om_1=\sum_mf^m_{1-\alpha}e^m_\alpha$ for each $\alpha\in(0,1)$. For the integration the point $\alpha=1$ does not count because $\CB{1}$ has measure $0$. By \it{dominated convergence}, we find
\bmun{
\int_0^1\AB{z_1x_\gamma,y_\gamma\om_1}\,d\gamma
~=~
\int_0^1\AB{z_1x_\gamma,y_\gamma\bfam{{\textstyle\sum_mf^m_{1-\gamma}e^m_\gamma}}}\,d\gamma
\\
~=~
\sum_m\int_0^1\AB{z_1x_\gamma,y_\gamma f^m_{1-\gamma}e^m_\gamma}\,d\gamma
~=~
\sum_m\int_0^1\AB{z_1,y_\gamma f^m_{1-\gamma}}\AB{x_\gamma,e^m_\gamma}\,d\gamma.
}\emun
Observe that $\norm{f^m_\alpha}\le1$ for all $\alpha\in(0,1),m\in\N$. There exists an $m$ such that $\int_0^1\snorm{f^m_{1-\gamma}}^2\,d\gamma<1$. (Indeed, if this integral is $1$ for a certain $m_0$, then it is $0$ for all other $m\ne m_0$.) Choose $X=e^m$ for that $m$, so that $\AB{x_\gamma,e^m_\gamma}=1$ for all $\gamma$,  and put $M:=\sqrt{\int_0^1\snorm{f^m_{1-\gamma}}^2\,d\gamma}$. Then
\bmun{
\Babs{\int_0^1\AB{z_1x_\gamma,y_\gamma\om_1}\,d\gamma}
~=~
\Babs{\int_0^1\AB{z_1,y_\gamma f^m_{1-\gamma}}\,d\gamma}
~=~
\Babs{\BAB{z_1,\int_0^1y_\gamma f^m_{1-\gamma}\,d\gamma}}
~\le~
\Bnorm{\int_0^1y_\gamma f^m_{1-\gamma}\,d\gamma}
\\
~\le~
\int_0^1\snorm{f^m_{1-\gamma}}\,\snorm{y_\gamma}\,d\gamma
~\le~
\sqrt{\int_0^1\snorm{f^m_{1-\gamma}}^2\,d\gamma}~\sqrt{\int_0^1\snorm{y_\gamma}^2\,d\gamma}
~=~
M\norm{Y}
~=~
M
~<~
1.
}\emun
The constant $M$ is independent of the choice of the unit vector $Y$. In conclusion, for $X=e^m$ we have $\norm{T_1(a)}\le M^2<1$. Therefore, $T_1(a)$ cannot be a non-zero projection. So, $T_1$ is not an endomorphism.\qed

\bobn
The Arveson system in the theorem is arbitrary. Since the $\vt$ constructed in \cite{Ske06} and the $\rho$ constructed in \cite{Flo08} coincide (for the choice of the parameters in $\rho$ as specified before) for integer $t=n\in\N_0$, also the compressed Markov maps $T_t$ coincide at least for integer $t=n$. As the theorem says $T_1$ is not an endomorphism, it follows that \cite[Theorem 4.4]{Flo08} is false. (We believe that the error is in Lemma 4.1. Check it for $t=s=1$, applying both sides to $\om_1x_1\otimes(Y\otimes\breve{z})$ when $x_1$ is taken from a unit $x^\otimes$ and $\om_1$ is taken from another unit $\om^\otimes$; see the computations below.)

For whom who wishes to have more concrete examples, we mention that it is possible to obtain simpler and calculable examples when the Arveson system $E^\otimes$ is spatial. In that case, we would choose a unital unit $\om^\otimes=\bfam{\om_t}_{t\in(0,\infty)}$ and for $\om_1$ really the member at $t=1$ of that unit. With this choice, the part $u_{n,\alpha}^*(y_\alpha\om_1^n)$ in $w_n^*(Y\otimes\breve{\om})$ may be computed as $y_\alpha\om_{n-\alpha}\otimes\om_\alpha$. When  computing $(\id_n\otimes X^*)\bfam{y_\alpha\om_{n-\alpha}\otimes\om_\alpha}_{\alpha\in\LO{0,1}}$ this gives $\int_0^1y_\alpha\om_{n-\alpha}\AB{x_\alpha,\om_\alpha}\,d\alpha$. Taking $x_\alpha=\om_\alpha$ and for $y_\alpha$ pieces from an independent unit, $\AB{Y,X}\AB{X,Y}$ and \eqref{D} can be computed. This works for an arbitrary spatial Arveson system of index not smaller than $1$. Of course, it also works for type I systems, that is, for Fock spaces (with $\om^\otimes$ the vacuum unit). Here, everything may be computed explicitly in terms of exponential vectors.
\eobn

\bremn
We should note that there is a simple theoretical argument, why a type III \nbd{E_0}semi\-group $\vt$ (that is, the Arveson system of $\vt$ is type III) can never be compressed to an automorphism semigroup $T$. (See the proof of the proposition below for the following terminology.) In fact, the Arveson system of $\vt$ contains the Arveson system of the minimal dilation of $T$, and the minimal dilation of an \nbd{E_0}semigroup (that is, in particular, of an automorphism semigroup) $T$ is $T$ itself. But the Arveson system of an automorphism semigroup would be the ``trivial'' one, $\bfam{\C}_{t\in(0,\infty)}$.%
\footnote[2]{\label{2}
Recall that, in these notes like Arveson in \cite{Arv89}, we did exclude the one-dimensional case. In fact, our Theorem is false in the one-dimensional case, and our proof breaks down once we have only one element in our measurable ONB.
}%
~ And the ``trivial'' Arveson system, like every Arveson system containing it, has a unit. This is not possible if $\vt$ is type III.

But our theorem is much more far-reaching. It tells that, no matter from which Arveson system $E^\otimes$ we start, $T$ is \hl{proper} in the sense that it is not even an endomorphism semigroup. In the remainder, we explain briefly why this promises to provide the first examples of nontrivial type III Markov semigroups.
\eremn

A \hl{type III} or \hl{nonspatial} Markov semigroup is a Markov semigroup with type III Arveson system. (This property is equivalent to the property that the semigroup has no \it{units} in the sense of Arveson \cite[Definition 2.1]{Arv97a}; see Bhat, Liebscher and Skeide \cite{BLS10}. It should not be confused with Powers' definition \cite{Pow04}, which is more restrictive.) Of course, every type III \nbd{E_0}semigroup is also an example for a type III Markov semigroup. By a \hl{nontrivial} type III Markov semigroup we understand a proper type III Markov semigroup that is not the tensor product of a type III \nbd{E_0}semigroup and a proper spatial Markov semigroup.

So far, there are no known examples of such nontrivial type III Markov semigroups. With some basic knowledge about minimal dilation and Arveson system of a Markov semigroup, our theorem allows to show that for certain type III Arveson systems, Floricel's Markov semigroup, necessarily type III, is nontrivial. The prerequisits are collected in the following proposition and its proof. Observe that with $Q$ also the projection $Q_t:=\vt_t(Q)$ is increasing for $\vt$. For $\alpha\ge0$, we denote by $T^\alpha$ the Markov semigroup on $Q_\alpha\sB(K)Q_\alpha$ obtained by compressing $\vt$ with $Q_\alpha$. Observe that with $T$, also $T^\alpha$ is proper. (This follows from $\vt_\alpha\circ T_t=T^\alpha_t\circ\vt_\alpha$. So, if $T_t$ does not factor on $a_1a_2$ ($a_i\in Q\sB(K)Q$), then $T^\alpha_t$ does not factor on $\vt_\alpha(a_1)\vt_\alpha(a_2)$ ($\vt_\alpha(a_i)\in Q_\alpha\sB(K)Q_\alpha$).)

\bpropn
Let $(\vt,Q)$ be a dilation of a Markov semigroup $T$.
\begin{enumerate}
\item
If $E^\otimes$ is an Arveson system that does not factor into the tensor product  of two Arveson systems, then, for each $\alpha>0$, $T^\alpha$ is a proper Markov semigroup that does not factor into the tensor product of two Markov semigroups.

\item
If $E^\otimes$ is an Arveson system that has no subsystem factoring into the tensor product  of two Arveson systems, then $T$ is a proper Markov semigroup that does not factor into the tensor product of Markov semigroups.
\end{enumerate}
\epropn

\proof
The dilation $(\vt,Q)$ of $T$ is \hl{minimal} if the smallest subspace of $K$ invariant for $\vt_t(a)$ $(t\in\R_+,a\in Q\sB(K)Q)$ and containing $QK$ is $K$. By Bhat \cite[Theorem 4.7]{Bha96}, every (normal) Markov semigroup $T$ on $\sB(H)$ admits a \it{minimal dilation} and that minimal dilation is unique up to suitable unitary equivalence. Bhat \cite[Section 6]{Bha96} defines the \hl{Arveson system of $T$} as the Arveson system of the unique minimal dilation. (This Arveson system can be constructed directly as explained in Skeide \cite{Ske03c} following the construction of Bhat and Skeide \cite{BhSk00}, or in Bhat and Mukherjee \cite{BhMu10} following notions of Arveson \cite{Arv97a}. But the statement we need here, really, is that the Arveson system of $T$ is that of the minimal dilation.)

There are two easy to verify consequences. Firstly, every dilation can be compressed to the smallest invariant subspace containing $QK$ (as above) to obtain the minimal dilation; see Bhat \cite[Section 3]{Bha01}. (See also Shalit and Solel \cite[Theorem 5.12]{ShaSo09} for a similar result in more general circumstances.) In either way to construct the Arveson system of that dilation (Arveson's \cite{Arv89} and Bhat's \cite{Bha01}) it is easy to see that the projection onto that subspace gives rise to a projection morphism of that Arveson system onto a subsystem that is the Arveson system of the minimal dilation: The Arveson system of every dilation contains the Arveson system of the minimal dilation. Secondly, given two Markov semigroups, the tensor product of their minimal dilations is the minimal dilation of their tensor product; this is mentioned in \cite{Bha96} between Theorems 6.3 and 6.4.

Putting these two statements together, immediately proves 2. (The Arveson system of $T$ is a subsystem of $E^\otimes$. If $E^\otimes$ has no subsystem that factors, then the minimal dilation does not factor, thus, neither does $T$.)

Statement 1 follows the same way from the following two theorems. \cite[Theorem 3.7]{Bha01}: If $(\vt,Q)$ is a \hl{primary dilation} (that is, if $Q_t\uparrow\id_K$), then, for all $\alpha>0$, the dilation $(\vt,Q_\alpha)$ of $T^\alpha$ is minimal. \cite[Theorem 3.6(ii)]{Bha01} (reformulated for our needs): If $(\vt,K)$ is not primary, then it has a corner containing $Q$ (hence, $Q_t$) that is a primary dilation with the same Arveson system as $\vt$, to which the former theorem can be applied. (Recall that, by Footnote \ref{2}, $\vt$ does not consist of automorphisms and \cite[Theorem 3.6(i)]{Bha01} does not apply. Anyway, without (the not very difficult direct) proof we communicate that the dilation $(\vt,Q)$ of $T$ as in our theorem, actually, is primary. For Floricel's dilation this statement is contained in \cite[Proposition 4.2]{Flo08}, and since $Q$ is increasing, it is sufficient to know it ony for integer times $t=n$, for which we clarified equality with \cite{Flo08}.)\qed

\bempn{Supplement.~}
If we specify that the Arveson (sub)system does not factor into certain types, then the Markov semigroup does not factor into these types either.
\eempn

\bcorn
If $E^\otimes$ is a type III Arveson system that does not factor into the tensor product of a type III system and a spatial system, then the semigroups $T^\alpha$ $(\alpha>0)$ derived from Floricel's dilation are nontrivial type III Markov semigroups. If $E^\otimes$ has even no subsystems factoring in that way, then Floricel's Markov semigroup itself is nontrivial type III.
\ecorn

Existence of such Arveson systems is, however, an open question. (Good candidates are \it{generalized CCR-flows} from Izumi and Srinivasan \cite{IzSr08} with one-dimensional multiplicity space.)

\lf\lf\noindent
\bf{Acknowledgments.~}
I would like to express my gratitude to Rajarama Bhat, Claus Köstler, and Volkmar Liebscher, with whom I had very useful discussions. I also wish to thank Rajarama Bhat for drawing my attention to \cite[Theorem 3.7]{Bha01}. Last but surely not least, I wish to thank the referee for careful reading and many valuable suggestions that improved the paper.


\vfill

\setlength{\baselineskip}{2.5ex}

\newcommand{\Swap}[2]{#2#1}\newcommand{\Sort}[1]{}
\providecommand{\bysame}{\leavevmode\hbox to3em{\hrulefill}\thinspace}
\providecommand{\MR}{\relax\ifhmode\unskip\space\fi MR }
\providecommand{\MRhref}[2]{%
  \href{http://www.ams.org/mathscinet-getitem?mr=#1}{#2}
}
\providecommand{\href}[2]{#2}

\lf\noindent
Michael Skeide: \it{Dipartimento S.E.G.e S., Università degli Studi del Molise, Via de Sanctis, 86100 Campobasso, Italy},
E-mail: \href{mailto:skeide@unimol.it}{\tt{skeide@unimol.it}},\\
Homepage: \url{http://www.math.tu-cottbus.de/INSTITUT/lswas/_skeide.html}

\end{document}